\documentclass[12pt, reqno]{amsart}
\usepackage{amsmath, amsthm, amscd, amsfonts, amssymb, graphicx, hyperref,color}
\newtheorem{theorem}{Theorem}[section]
\newtheorem{lemma}[theorem]{Lemma}
\newtheorem{proposition}[theorem]{Proposition}
\newtheorem{corollary}[theorem]{Corollary}
\theoremstyle{definition}
\newtheorem{definition}[theorem]{Definition}
\newtheorem{example}[theorem]{Example}

\theoremstyle{remark}
\newtheorem{remark}[theorem]{Remark}
\numberwithin{equation}{section}

\begin{document}
\setcounter{page}{1}


\title[The Moore-Penrose inverse of accretive  operators ]
      {The Moore-Penrose inverse of accretive  operators with  application to quadratic operator pencils }
      
\author[ F. Bouchelaghem, M. Benharrat]{Fairouz Bouchelaghem$^1$, Mohammed Benharrat$^{2*}$}
\address{$^{1}$ 
D\'{e}partement de Math\'{e}matiques, Univérsit\'{e} Oran 1-Ahmed Ben Bella, BP 1524 Oran-El M'naouar, 31000 Oran, Alg\'{e}rie.}
\email{\textcolor[rgb]{0.00,0.00,0.84}{fairouzbouchelaghem@yahoo.fr}}
\address{$^{2}$ D\'{e}partement de G\'{e}nie des syst\'{e}mes,
	Ecole Nationale Polytechnique d'Oran-Maurice Audin (Ex. ENSET d'Oran), 
	BP 1523 Oran-El M'naouar, 31000 Oran, Alg\'{e}rie.}
\email{\textcolor[rgb]{0.00,0.00,0.84}{	mohammed.benharrat@enp-oran.dz, mohammed.benharrat@gmail.com}}





\date{15/01/2021.
	\newline \indent $^{*}$ Corresponding author\\
	This work was supported by the Laboratory of Fundamental and Applicable Mathematics of Oran (LMFAO) and the Algerian research project: PRFU, no. C00L03ES310120180002.}
\begin{abstract}
We establish some relationships between an m-accretive operator and its Moore-Penorse inverse. We derive some perturbation result of the Moore-Penorse inverse of a maximal accretive operator. As an application  we give a  factorization theorem for  a quadratic  pencil of accretive operators. Also, we study a result of existence, uniqueness, and maximal regularity of the strict solution  for complete abstract second order differential equation. Illustrative examples are also given. 
\end{abstract}
 \maketitle

\section{introduction}
The  Moore-Penrose inverse of a linear operator in Hilbert space is a useful generalization of the
ordinary inverse. This generalized  inverse is an important theoretical and practical tool in algebra and
analysis (Markov chains, singular differential and difference equations, iterative methods...), see \cite{israel,Wang2018}.  In particular, in \cite{Huang2011,Huang2012} (and the references therein) the perturbation analysis  for the Moore-Penrose inverse of closed  operators has been considered. Also,  the expression of the generalized inverse of the perturbed operator has been investigated. In the paper  \cite{Kurmayya2016} necessary and sufficient conditions for the cone
nonnegativity of Moore–Penrose inverses of unbounded Gram operators are
derived. These conditions include statements on acuteness of certain closed
convex cones in infinite-dimensional real Hilbert spaces. In \cite{Benharrat2020} a complete description of the left quotient and the right quotient of two bounded operators operators is given via the Moore-Penrose inverse.  The objectives of this
paper are to derive the properties of m-accretive operators via the Moore-Penrose inverse and establish some interesting results, especially for the perturbation analysis of Moore-Penrose inverses as well as  of maximal accretive operators. Recall that  a linear operator $T$ with domain $\mathcal{D}(T)$ in a complex Hilbert space $\mathcal{H}$, is called accretive if its numerical range $W(T)$ is contained in the closed right half-plane, and if further has no proper accretive extensions in $\mathcal{H}$, it called maximal accretive, m-accretive for short. In particular, every m-accretive operator is accretive
and  closed densely defined, its adjoint is also m-accretive (cf. \cite{Kato}, p. 279).  This class is of particular interest and related to the  semi-group theory in follwing sens: an operator $T$ is m-accretive if and only if $-T$ generates a strongly continuous contraction semigroup (Theorem  of Lumer-Phillips).

In this paper, we explore the following two questions, the first what can be said about the m-accretivity of the Moore-Penrose of an m-accretive operator and conversely? the second concern the perturbation problem: Let $T$ be m-accretive operator with a bounded Moore-Penrose inverse, what condition on
the  operator $S$  can guarantee that  $T+S$ is m-accretive and its Moore-Penrose inverse 
exists and it has the simplest expression? 

In this work, we give a certain answers to the mentioned problems.
This paper is organized as follows: In section 2,
we establish some relationships between an m-accretive operator and its Moore-Penorse inverse. In section 3, we consider the perturbation for the  m-accretive operator and  its Moore-Penrose inverse. We prove that under weaker conditions that considered perturbation does not change
the null space and the range space, consequently the perturbed operator is  a closed EP operator. Utilizing this result, we study a class of a
quadratic operator pencil $\mathcal{Q} (\lambda)=\lambda^2 I-2 \lambda T-S,$ $(\lambda \in \mathbb{C})$, where the coefficients of which are  accretive  operators. Our aim, in Section 4, is   to investigate   a canonical factorization like $(\lambda I- Z_1)(\lambda I-Z_2)$  for of such pencils based on the  perturbation theory of  accretive operators. We also obtain a
criterion in order that the linear factors, into which
the pencil splits,  generates an holomorphic semi-group
of contraction operators. As an illustration, in section 5, we  establish a theorem  of existence, uniqueness, and maximal regularity of the strict solution  of an abstract second order evolutionary equations generated by such pencils. 
 




\section{Accretive  operator  and the Moore-Penrose inverse}
Throughout this paper  $\mathcal{H}$ is a complex   Hilbert space with inner product $<\cdot, \cdot>$ and norm $\|\cdot\|$. For a closed linear operator $T$ on $\mathcal{H}$ we denote by $\mathcal{D} (T)$, $\mathcal{R}(T)$, $\mathcal{N}(T)$, $\sigma (T)$ and $\rho (T)$
the domain, the range, the kernel, the spectrum and the resolvent set of $T$, respectively.  The space of
bounded linear operators on $\mathcal{H}$ is denoted by $\mathcal{B} (\mathcal{H})$. For two possibly unbounded
linear operators $T$, $S$ on $\mathcal{H}$ their product $TS$ is defined on its natural domain
$\mathcal{D} (TS) := \{x \in  \mathcal{D} (S) :  Sx \in \mathcal{D} (T)\} $ and their sum $T+S$ is defined in $\mathcal{D} (T+S)  = \mathcal{D} (T) \cap \mathcal{D} (S)$. An inclusion $T \subseteq S$ denotes inclusion of graphs, i.e., it means that $S$ extends $T$. A possibly unbounded operator
$T$ on $\mathcal{H}$ commutes with a bounded operator $S\in \mathcal{B} (\mathcal{H})$ if the graph of $T$ is $S \times S$-invariant,
or equivalently if $ST \subseteq TS$.


Recall that  a linear operator $T$ with domain $\mathcal{D}(T)$ in  $\mathcal{H}$ is said to be accretive if
$$
{\rm Re}<Tx,x>\geq 0 \qquad \text{ for all } x\in \mathcal{D}(T)
$$
or, equivalently if
$$\Vert(\lambda +T)x\Vert\geq \lambda \Vert x\Vert \qquad \text{ for all } x\in \mathcal{D}(T) \text{ and  }\lambda>0.$$
An accretive operator $T$ is called \textit{maximal accretive}, or
\textit{$m$-accretive} for short, if  $T$ has no proper accretive extensions in $\mathcal{H}$. The following conditions are equivalent:
\begin{enumerate}
	\item $T$ is $m$-accretive.
	\item  $$ (\lambda + T)^{-1} \in \mathcal{B} \ (\mathcal{H}) \quad  \text{ and } \quad \left\|(\lambda + T)^{-1}\right\| \leq \frac{1}{\lambda}   \text{ for } \lambda >0.$$
	\item $T$ is accretive densely defined and $\mathcal{R}(\lambda +T)=\mathcal{H}$ for some (and hence for every) $\lambda>0$;
	\item $T$ is accretive densely defined and closed, and $T^*$ is accretive;
	\item $-T$ generates contractive one-parameter semigroup $\mathcal{T}(t)=\exp(-tT)$, $t\ge 0$.
\end{enumerate}
In particular, a bounded accretive  operator  is m-accretive. 

The numerical range of a linear operator  $T: \mathcal{D}  (T)\rightarrow  \mathcal{H}$  it is defined by
\begin{equation}
W(T):=\{ <Tx, x> : \quad  x\in \mathcal{D}  (T) , \quad \text{ with } \left\|x\right\|=1\},
\end{equation}
It is well-known that  $W(T)$ is a convex set of the complex plane (the Toeplitz-Hausdorff
theorem), and in general is neither open nor closed, even for a closed operator $T$. Clearly,  an operator $T$  is accretive when   $W(T)$ is contained in the closed right half-plane
\[
W(T) \subset \overline{\mathbb{C}_+} := \{ z\in\mathbb{C} :  {\rm Re} (z) \ge 0\}.
\]
Further,  if $T$ is m-accretive operator then $W(T)$ has the so-called spectral inclusion property
\begin{equation}\label{specinclusion}
\sigma (T)\subset\overline{ W(T)}.
\end{equation}

Recall that a linear operator $T$ in  $\mathcal{H}$ is called \emph{sectorial} with vertex $z=0$ and
semi-angle $\omega \in [0,\pi/2)$, or \emph{$\omega$-accretive} for short, if its numerical range is contained in a closed sector with semi-angle~$\omega$,
\begin{equation}
\label{sectorial}
W(T) \subset \overline{\mathcal{S}(\omega)}:=\left\{z\in\mathbb{C} : |\arg z|\leq \omega \right\}
\end{equation}
or, equivalently, $$|Im<Tx, x> |\!\le\! \tan\omega \, {\rm Re}<Tx, x> \qquad \text{ for all }  x \in  \mathcal{D}  (T).$$
An $\omega$-accretive  operator  $T$ is called   m-$\omega$-accretive, if it is $m$-accretive. We have $T$ is m-$\omega$-accretive if and only if the operators $e^{\pm i\theta} T$ is  m-accretive for $\theta=\frac{\pi}{2}-\omega$, $0 < \omega\leq  \pi/2$.
The resolvent set of an m-$\omega$-accretive operator $T$ contains the
set $\mathbb{C} \setminus \overline{\mathcal{S}(\omega)}$ and
\[
\|(T-\lambda I)^{-1}\| \le \cfrac{1}{{\rm dist}\left(\lambda,\mathcal{S}(\omega)\right)},
\quad \lambda\in \mathbb{C}\setminus \overline{\mathcal{S}(\omega)}.
\]
In particular,  m-$\pi/2$-accretivity means m-accretivity.  A $0$-accretive operator
is symmetric. An operator is positive if and only if it is m-$0$-accretive. 

It is known that the $C_0$-semigroup $\mathcal{T}(t)=\exp(-tT)$, $t\ge 0$, has contractive  and holomorphic continuation into the sector $ \overline{\mathcal{S}(\pi/2\omega)}$ if and only if the generator $T$  is m-$\omega$-accretive, see \cite[Theorem V-3.35] {Kato}.

Recall that for bounded operator  $T$, we have
$$ 
{\rm Re}(T)=\dfrac{1}{2}(T+T^*) \qquad \text{ and } \qquad {\rm Im(T)}=\dfrac{1}{2i}(T-T^*), 
$$
where ${\rm Re}(T)$ and ${\rm Im(T)}$ are self-adjoint operators and called it the real and imaginary parts of $T$, with
$$
T={\rm Re}(T)+i{\rm Im(T)},
$$ Such  decomposition is unique and  called the cartesian decomposition of $T$. In this case, $T$ is accretive if    ${\rm Re}(T)$ is a nonnegative operator.

The spectral radius and the numerical radius of a bounded operator $T$  are  defined, respectively,  by
$$r(T)=  \sup_{\lambda \in \sigma (T)} \left| \lambda \right| $$
and
$$w(T)=\sup_{\left\|x \right\| =1} \left| <Tx,x>\right| .$$

The next results is a generalization  of \cite[Lemma 2. and Theorem 2.]{Arlin2003} from matrices to bounded operators.

\begin{lemma}\label{lambda0} Let $T\in \mathcal{B}(\mathcal{H})$ such that ${\rm Re}(T)\geq  \delta I$, for some  $\delta >0$. Denote by $\mathcal{S}={\rm Im(T)}({\rm Re}(T))^{-1}$. Then there exists $\lambda_0 \in \sigma (\mathcal{S}) $, $\left| \lambda_0\right|=r(\mathcal{S})=w(\mathcal{S})=\left\| \mathcal{S}\right\|$  such that   $T$ is a sectorial operator  with 
	semiangle $\omega=arctan\left| \lambda_0\right|$ and
	\begin{equation}\label{lambda00}
	\left| \lambda_0\right|\leq\sqrt{ \dfrac{\left\| T\right\|^2}{ \delta^2}-1}.
	\end{equation}
\end{lemma}
\begin{proof}By assumption,
	$$
	<{\rm Re}(T)x,x> \geq \delta \left\| x\right\|^2 
	$$
	for all $x\in \mathcal{H}$. Hence we get
	$$
	\left|	<{\rm Im(T)}x,x> \right| \leq \left\| {\rm Im(T)}\right\| \left\| x\right\|^2\leq \dfrac{1}{\delta}\left\| {\rm Im(T)}\right\| <{\rm Re}(T)x,x>.
	$$	
	Thus, $T$ is sectorial operator with a vertex at the origin with a
	semiangle $\omega$. In particular, for ${\rm Re}(T)^{-\frac{1}{2}}x$, we have
	$$
	\left|	<{\rm Im(T)}{\rm Re}(T)^{-\frac{1}{2}}x,{\rm Re}(T)^{-\frac{1}{2}}x> \right| \leq \left\| {\rm Re}(T)^{-\frac{1}{2}}{\rm Im(T)}{\rm Re}(T)^{-\frac{1}{2}}\right\| \left\| x\right\|^2
	$$
	this implies that $tan(\omega)=\left\| \mathcal{T}\right\|$, where $\mathcal{T}={\rm Re}(T)^{-\frac{1}{2}}{\rm Im(T)}{\rm Re}(T)^{-\frac{1}{2}}$. 
	Since $\mathcal{T}$ is a self-adjoint operator we have $r(\mathcal{T})=w(\mathcal{T})=\left\| \mathcal{T}\right\|$. Thus  assert the existence of a $\lambda_0 \in \sigma(\mathcal{T}) $ such that $ \left| \lambda_0\right|=r(\mathcal{T})=w(\mathcal{T})=\left\| \mathcal{T}\right\|$. Since
	$$\mathcal{T}={\rm Re}(T)^{-\frac{1}{2}}\mathcal{S} {\rm Re}(T)^{\frac{1}{2}}={\rm Re}(T)^{\frac{1}{2}}\mathcal{S}^* {\rm Re}(T)^{-\frac{1}{2}},$$
	the self-adjoint  operators $\mathcal{T}$ and $\mathcal{S}$ have the same spectrum  (which is real ) and hence have the same closure of the numerical range. This shows that $\lambda_0 \in \sigma (\mathcal{S}) $, with $\left| \lambda_0\right|=w(\mathcal{T})=w(\mathcal{S}) $ and hence $\omega=arctan\left| \lambda_0\right|$. 
	
	Now, assume that $\lambda_0$ is an eigenvalue of $\mathcal{T}$ with $ \left| \lambda_0\right|=\left\| \mathcal{T}\right\|$, there exists $u\in \mathcal{H}$ with $ \left\| u\right\|=1$ and $\mathcal{T}u=\lambda_0 u $. We have
	 $$  {\rm Re}(T)^{-\frac{1}{2}}T{\rm Re}(T)^{-\frac{1}{2}}=  (I+i \mathcal{T} )$$
	and
$$
	\left\| {\rm Re}(T)^{-\frac{1}{2}}T{\rm Re}(T)^{-\frac{1}{2}}\right\|^2\geq \left\| (I+i \mathcal{T} )u\right\|^2
	=1+ \left\| \mathcal{T}\right\|^2=1+\left| \lambda_0\right|^2,
$$
which implies that 
\begin{equation*}
\left| \lambda_0\right|^2\leq \dfrac{\left\| T\right\|^2}{ \delta^2}-1.
\end{equation*}
Now we consider the general case. Let $\varepsilon >0$. It follows from the spectral theorem that there exists a self-adjoint bounded operator $P$ such that $ \left\| P\right\|\leq \varepsilon$ and the operator $\mathcal{T}+P $ has an eigenvalue such that the modulus equals $ \left\|\mathcal{T}+ P\right\|$. As above, we take  the operator $I + i(\mathcal{T}+P)$ instead of $(I+i \mathcal{T} )$, we get
$$ \dfrac{\left\| T\right\|^2}{ \delta^2} +\varepsilon^2\geq 1+\left| \lambda_0\right|^2.$$
Letting $\varepsilon \longrightarrow 0$, we obtain \eqref{lambda00}.
\end{proof}

If we assume, in  Lemma \ref{lambda0}, that the numerical range of $T$ is closed, then $\lambda_0\in W(T)$, but the extreme points of the numerical
range are in the point spectrum,  so $\lambda_0$ must be an eigenvalue of $\mathcal{T}$. So, we have the following 
\begin{corollary}\label{lambda0pp} Let $T\in \mathcal{B}(\mathcal{H})$ with closed numerical range such that ${\rm Re}(T)\geq  \delta I$, for some  $\delta >0$. Denote by $\mathcal{S}={\rm Im(T)}({\rm Re}(T))^{-1}$. Then there exists $\lambda_0 \in \sigma_p (\mathcal{S}) $, $\left| \lambda_0\right|=w(\mathcal{S}) $ such that   $T$ is a sectorial operator  with 
	semiangle $\omega=arctan\left| \lambda_0\right|$ and $\lambda_0$ verifies \eqref{lambda00}.
\end{corollary}
\begin{remark}
	\begin{enumerate}
		\item Since ${\rm Re}(T)$ is strongly nonnegative, we know that 
		$$ 
		\sigma (\mathcal{S})\subseteq \{ \dfrac{\beta}{\alpha} : \quad \alpha \in \overline{W({\rm Re}(T))},  \beta \in \overline{W({\rm Im(T)})}\}.
		$$
		So $ \lambda_0=\dfrac{\beta_0}{\alpha_0}$ for  $\alpha_0 \in \overline{W({\rm Re}(T))}$ and   $\beta_0 \in \overline{W({\rm Im(T)})}$.
		\item As we can see from the  proof above that $T$ can be represented as
		$$ T=  {\rm Re}(T)^{\frac{1}{2}} (I+i \mathcal{T} ){\rm Re}(T)^{\frac{1}{2}}$$
		with $\mathcal{T}={\rm Re}(T)^{-\frac{1}{2}}{\rm Im(T)} {\rm Re}(T)^{-\frac{1}{2}}$  and $\tan (\omega)=\left\| \mathcal{T}\right\|$. This is exactly   the representation given in \cite[Theorem VI-3.2]{Kato}. In our case  the selfadjoint operator is uniquely determined.
	\end{enumerate}
	
\end{remark}

Next, in order to  give some new results  about accretive  operator by using the Moore-Penrose inverse, let recall  the definition of this generalized inverse for a closed densely defined operator.
\begin{definition}\cite{israel}\label{MPInv}
	Let $T$ be  a closed densely defined on $\mathcal{H}$. Then there exists a unique closed densely defined
	operator $T^{\dag}$,   with domain $\mathcal{D}(T^{\dag})= \mathcal{R}(T) \oplus \mathcal{R}(T)^\bot$  such that 
	$$ TT^{\dag}T=T \quad \text {on }  \mathcal{D}(T), \qquad T^{\dag}TT^{\dag}=T^{\dag} \quad \text {on }  \mathcal{D}(T^{\dag}), $$
	$$ TT^{\dag}=P_{\overline{\mathcal{R}(T)}} \quad \text { on }  \mathcal{D}(T^{\dag}), \qquad T^{\dag}T=P_{\mathcal{N}(T)^{\bot}} \quad \text { on }  \mathcal{D}(T), $$
	with $P_{\mathcal{M}}$ denotes the orthogonal projection onto a closed subspace $\mathcal{M}$.
\end{definition}
This unique operator $T^{\dag}$ is called the Moore-Penrose inverse of $T$.  (or  the Maximal
Tseng generalized Inverse in the terminology of \cite{israel}). Clearly, 
\begin{enumerate}
	\item $\mathcal{N}(T^{\dag})=  \mathcal{R}(T)^\bot$,
	\item  $\mathcal{R}(T^{\dag})=  \mathcal{N}(T)^\bot\cap \mathcal{D}(T)$.
\end{enumerate}
As a consequence of the closed graph theorem  $T^{\dag}$ is bounded if and only if  $\mathcal{R}(T)$ is closed in $\mathcal{H} $, see \cite{israel}.

Now, if we assume that $T$ is an m-accretive operator, then
\begin{equation}
\label{kernel-accretive}
\mathcal{N} (T) = \mathcal{N} (T^*) \quad  \text{ and } \quad \mathcal{N} (T)\subseteq\mathcal{D} (T) \cap \mathcal{D} (T^*).
\end{equation} 
Thus   $\overline{\mathcal{R}(T)}=\overline{\mathcal{R}(T^*)}$ and $\mathcal{H}=\overline{\mathcal{R}(T)}\oplus \mathcal{N}(T)$. Consequently,  the operator  $T$ is written in a matrix form with respect to mutually orthogonal subspaces decomposition as follows
$$T =\begin{bmatrix}
T_1& 0 \\
0& 0
\end{bmatrix} : \begin{bmatrix}
\overline{\mathcal{R}(T)}\\
\mathcal{N}(T)
\end{bmatrix}\longrightarrow  \begin{bmatrix}
\overline{\mathcal{R}(T)}\\
\mathcal{N}(T)
\end{bmatrix};
$$
with $T_1$ is an  operator on  $\overline{\mathcal{R}(T)} \cap \mathcal{D}(T)$ is injective with dense range in $\overline{\mathcal{R}(T)}$. Also,  its Moore-Penrose inverse is given by
$$  T^{\dag} =\begin{bmatrix}
T_1^{-1}& 0\\
0 & 0
\end{bmatrix}: \begin{bmatrix}
\mathcal{R}(T)\\
\mathcal{N}(T)
\end{bmatrix}\longrightarrow  \begin{bmatrix}
\overline{\mathcal{R}(T)}\\
\mathcal{N}(T)
\end{bmatrix},$$
with $T_1^{-1}$ from $\mathcal{R}(T)$ to $\overline{\mathcal{R}(T)}\cap \mathcal{D}(T)$  is closed operator densely defined  on  $\overline{\mathcal{R}(T)}$ and $\mathcal{N}(T^{\dag})=\mathcal{N}(T)=\mathcal{N}(T^*)$. Further,  $\mathcal{R}(T)$ is closed if and only if  $T_1^{-1}$ is bounded from $\mathcal{R}(T)$ to $\mathcal{R}(T)\cap \mathcal{D}(T)$.

In the next result, we redefine  a bounded sectorial operator  with 
semiangle $<\pi/2$ via the Moore-Penrose inverse.

\begin{lemma}\label{sec1} Let $T\in \mathcal{B}(\mathcal{H})$.    $T$ is a sectorial operator  with 
	semiangle $\omega$, $0 \leq \omega<  \pi/2$,  if and only if the following two conditions are fulfilled:
	\begin{enumerate}
		\item ${\rm Re}(T)\geq 0$.
		\item $\overline{\mathcal{R}(T)}\subset \overline{\mathcal{R}({\rm Re}(T))}$.
	\end{enumerate}
	In this case, $\omega=arctan\left| \lambda_0\right|$ for some $\lambda_0 \in \sigma (\mathcal{S}) $, $\left| \lambda_0\right|=w(\mathcal{S})$ where  $\mathcal{S}={\rm Im(T)}({\rm Re}(T))^{\dag}$.
\end{lemma}
\begin{proof}Let $T\in \mathcal{B}(\mathcal{H})$ be  a sectorial operator  with semiangle $\omega$. Therefore ${\rm Re}(T)$ is a nonnegative and $\mathcal{N}(T)=\mathcal{N}(T^*)=\mathcal{N}({\rm Re}(T))$. It follows that  $\overline{\mathcal{R}(T)}= \overline{\mathcal{R}({\rm Re}(T))}$.
	
	Conversely, let the two conditions of the theorem hold. Then $\mathcal{N}(T)=\mathcal{N}(T^*)\subset \mathcal{N}({\rm Re}(T))$, indeed, if $Tx=0$ for some $x\neq 0$, then ${\rm Re}<Tx,x>=<{\rm Re}(T)x,x>=0$. Consequently, $ {\rm Re}(T)x=0$ and $T^*x=2{\rm Re}(T)x-Tx=0$. Since  $\mathcal{H}=\mathcal{N}({\rm Re}(T)) \oplus \overline{\mathcal{R}({\rm Re}(T))}$,  we get 
	$$
	\overline{\mathcal{R}({\rm Re}(T))}\subset \overline{\mathcal{R}(T)}= \overline{\mathcal{R}(T^*)}
	$$
	with  $ \overline{\mathcal{R}(T)}$ reduces $T$. By (2) we obtain $\overline{\mathcal{R}({\rm Re}(T))}= \overline{\mathcal{R}(T)}= \overline{\mathcal{R}(T^*)}.$ It follows that the subspace $\overline{\mathcal{R}(T)}$ reduces also the operator ${\rm Re}(T)$. Moreover, the restriction of ${\rm Re}(T)$ to  $\overline{\mathcal{R}(T)}$  is strongly nonnegative. So by Lemma \ref{lambda0}, the  restriction of $T$ to  $\overline{\mathcal{R}(T)}$ is sectorial operator with semiangle  $\omega=arctan\left| \lambda_0\right|$ such that $\lambda_0 \in \sigma (\mathcal{S}_{| \overline{\mathcal{R}(T)}}) $ and $\left| \lambda_0\right|=\sup_{\left\|x \right\| =1} \left| <\mathcal{S}_{| \overline{\mathcal{R}(T)}}x,x>\right| $, where
	$\mathcal{S}_{| \overline{\mathcal{R}(T)}}={\rm Im(T)}({\rm Re}(T)_{| \overline{\mathcal{R}(T)}})^{-1}P_{| \overline{\mathcal{R}(T)}}$. In this case $\mathcal{S}={\rm Im(T)}({\rm Re}(T))^{\dag}$.
	
\end{proof}
Now, we consider unbounded operator $T$. 
\begin{proposition} If $T$ is  m-accretive operator, then $T^{\dag}$ is m-accretive.
\end{proposition}
\begin{proof}By assumption,
	$${\rm Re}<Tx,x>\geq 0 \quad  \text{ for all } x\in  \mathcal{D}(T)\cap \mathcal{N}(T)^\bot=\mathcal{R}(T^{\dag}). $$
	 Hence
	$${\rm Re}<y, T^{\dag}y>\geq 0 \quad  \text{ for all } y\in  \mathcal{R}(T). $$
	 Now let $x\in \mathcal{D}(T^{\dag})= \mathcal{R}(T) \oplus \mathcal{N}(T)$, then $x=x_1+x_2$, with $x_1\in \mathcal{R}(T)$ and $x_2\in \mathcal{N}(T) $. Therefore,
$${\rm Re}<x,T^{\dag}x >={\rm Re}<x_1,T^{\dag}x_1 >\geq  0, $$
which implies
$${\rm Re}<x,T^{\dag}x >\geq 0 \quad  \text{ for all } x\in \mathcal{D}(T^{\dag}).$$
Since $T^{\dag}$ is closed densely defined and $(T^{\dag})^*$ is accretive, it follows that $T^{\dag}$ is m-accretive.
\end{proof}
It well known that by \cite[Theorem 2; p 341]{israel},  $T^{\dag \dag}=T$, this yields to
\begin{corollary}  $T^{\dag}$ is  m-accretive operator if and only if  $T$ is m-accretive.
\end{corollary}
\begin{corollary}   $T$ is  m-accretive operator with closed range if and only if  $T^{\dag}$ is  bounded and accretive.
\end{corollary}
\begin{corollary}\label{EPAcc}   If $T$ is  m-accretive operator with closed range,  then $T$ is an EP (Equal Projections) operator, that is,  $T^{\dag}$ bounded and $TT^{\dag}=T^{\dag}T$ on $\mathcal{D}(T)$.
\end{corollary}
\begin{proposition}\label{thm:unitary} Let $T$ be an accretive bounded operator. If $W(T)\subseteq \overline{\mathbb{D}}$    and $W(T^\dagger)\subseteq \overline{\mathbb{D}}$ , then  $T$ is unitary on $\mathcal{R}(T)$.
\end{proposition}
\begin{proof} 
	It well known by \cite[Theorem 1.3-1]{GustafsonR97} that the numerical radius is equivalent to the usual operator norm;
	$$w(T)\leq \left\|T \right\|\leq 2w(T).$$
	Hence, the assumption that $w(T^\dagger)\leq 1$ implies that $T^\dagger$ is bounded. Thus  $\mathcal{R}(T)$ is closed. Since $T$ is m-accretive, then  $\mathcal{R}(T)=\mathcal{R}(T^*)=\mathcal{N}(T)^\bot$. We consider the restriction of $T$ from $\mathcal{R}(T)$ into itself. Since $T^\dagger_{|\mathcal{R}(T)}= (T_{|\mathcal{R}(T)})^{-1} $, $w((T_{\mathcal{R}(T)})^{-1})=w(T^\dagger_{|\mathcal{R}(T)})=w(T^\dagger)\leq 1$. Combining this with $w(T_{|\mathcal{R}(T)})=w(T)\leq 1$ and applying \cite[Corollary 1.]{Stampfli67} to $(T_{|\mathcal{R}(T)})^{-1}$ and $T_{|\mathcal{R}(T)}$, we conclude that $T_{|\mathcal{R}(T)}$ is unitary on $\mathcal{R}(T)$.
\end{proof}
\section{A perturbation results}

In the following, we shall consider the perturbation of  Moore-Penrose inverse of m-accretive operators.  This gives an example of EP-operators, see \cite[Theorem 3.12]{Qianglian2012}. The part about the Moore-Penrose invertibility of the operator $T+S$ appears to be new.
\begin{theorem}\label{MPInvsum}
	Let $T$ is  m-accretive operator and $S$ is bounded and accretive. We have
	\begin{enumerate}
		\item $T+S$ is m-accretive.
		\item If $\mathcal{R}(T)$ is closed, $\mathcal{R}(S)\subseteq \mathcal{R}(T)$ and $\left\| T^\dagger S\right\|  <1$.
	Then \begin{itemize}
		\item $\mathcal{R}(T+S)=\mathcal{R}(T)$ is closed and $\mathcal{N}(T+S)=\mathcal{N}(T)$.
		\item $T+S$  is an  EP  operator, and
		$$(T+S)^\dagger=(I+T^\dagger S)^{-1}T^\dagger =T^\dagger
		(I+ST^\dagger)^{-1}.$$ 
		In particular,  $$T^\dagger =(T+S)^\dagger
		(I+ST^\dagger),$$
	and
	 $${\left\| (T+S)^\dagger -T^\dagger \right\| \leq \displaystyle \frac{\left\| S\right\| \left\| T^\dagger\right\| ^2}{1-\left\| T^\dagger S\right\| }}.$$
	 \item Further, if  $S$ is $\theta$-accretive operator with $\theta\in[0,\frac{\pi}{2})$, then
	 \[
	 \left\| (T+S)^{\dagger}\right\|  \leq 2  \left\| T^{\dagger}\right\|  + (1+\tan\theta)^2  \, \left\| T^{\dagger}\right\|^2.
	 \]
	\end{itemize}	
	\end{enumerate}
	
\end{theorem}
\begin{proof} (1) Clearly, the operator $T+S$,  with $\mathcal{D}(T+S)=\mathcal{D}(T)$, is densely defined, closed and accretive. Since also its adjoint operator $(T+S)^*=T^* +S^*$ is accretive, the operator $T+S$ is m-accretive.

	(2) If $\mathcal{R}(S)\subseteq \mathcal{R}(T)$,  then it is obvious that $\mathcal{R}(T+S)\subseteq \mathcal{R}(T)$ and  $TT^\dagger S =P_{\mathcal{R}(T)}S=S$. Conversely, let $y\in \mathcal{R}(T)$, so $y=Tx$ for some $x\in \mathcal{D}(T)$. The condition $\left\| T^\dagger S\right\|  <1$ implies that $(I+T^\dagger S)^{-1}$ exists and bounded.
	Hence, there exists  $u\in \mathcal{D}(T)$ such that  $x=(I+T^\dagger S)u$. This  shows that  $y=T(I+T^\dagger S)u=Tu+Su \in \mathcal{R}(T+S)$. Hence $R(T)\subseteq \mathcal{R}(T+S)$. Consequently, $ \mathcal{R}(T+S)=\mathcal{R}(T)$ is closed.
	
	Since $T$ and $T+S$  are m-accretive with closed ranges, then	$$\mathcal{N}(T+S)=\mathcal{R}(T+S)^{\bot} =\mathcal{R}(T)^{\bot}= \mathcal{N}(T). $$
	
	Now we  prove that  $(T+S)^\dagger=(I+T^\dagger S)^{-1}T^\dagger$. Since,   $\mathcal{R}(T+S)$ is closed and  $\mathcal{N}(T+S)=\mathcal{R}(T+S)$ , by Corollary \ref{EPAcc},  it follows that $T+S$  is an  EP  operator.
	
	 Put $\mathbf{T}=(I+T^\dagger S)^{-1}T^\dagger$. We show that $\mathbf{T}$ satisfies all
	the axioms of the Definition \ref{MPInv}. First let us remark that, since $(I+T^\dagger S)^{-1}$ is invertible, $\mathcal{D}(\mathbf{T})=\mathcal{D}(T^\dagger)=\mathcal{R}(T) \oplus \mathcal{N}(T)$, $\mathcal{N}(\mathbf{T})=\mathcal{N}(T^{\dagger})=\mathcal{R}(T)^{\bot}=\mathcal{N}(T+S)$. 
	
	Let $v\in \mathcal{R}(\mathbf{T})$, then there exists  $u\in \mathcal{D}(\mathbf{T})$ such
	that $v=\mathbf{T}u=(I+T^\dagger S)^{-1}T^\dagger u$. Hence $T^\dagger
	u=v+T^\dagger Sv\in \mathcal{R}(T)\cap  \mathcal{D}(T)$. So $v=T^\dagger u-T^\dagger Sv \in \mathcal{D}(T)$.
	
	Now for $v\in \mathcal{D}(T)$,
	\begin{align*}
	\mathbf{T}(T+S)v&=(I+T^\dagger S)^{-1}T^\dagger (T+S)v\\
	&=(I+T^\dagger S)^{-1}T^\dagger (T+TT^\dagger S)v \quad (\text{since} \; S=TT^\dagger S)\\
	&=(I+T^\dagger S)^{-1}T^\dagger T(I+T^\dagger S )v\\
	&=(I+T^\dagger S)^{-1}P_{\mathcal{N}(T)^{\bot}}(I+T^\dagger S )v\\
    &=(I+T^\dagger S)^{-1}P_{\mathcal{R}(T)}(I+T^\dagger S )v\\
	&=(I+T^\dagger S)^{-1}(I+T^\dagger S )v\\
	&=v=P_{\mathcal{R}(T)}v\\
	&=P_{\mathcal{R}(T+S)}v\\
	&=P_{\mathcal{N}(T+S)^{\bot}}v.
	\end{align*}
and for $u\in \mathcal{D}(\mathbf{T})$,
	\begin{align*}
	(T+S)\mathbf{T}u=&(T+S)(I+T^\dagger S)^{-1}T^\dagger u \\
	=&(T+TT^\dagger S)(I+T^\dagger S)^{-1}T^\dagger u  \quad (\text{since} \; S=TT^\dagger S)\\
	=&T(I+T^\dagger S)(I+T^\dagger S)^{-1}T^\dagger u \\
	=&TT^\dagger u=P_{\mathcal{R}(T)}u\\
	=&P_{\mathcal{R}(T+S)}u.
	\end{align*}
	
	The uniqueness of $(T+S)^\dagger$ follows  from Definition \ref{MPInv}.
	
	Since $\mathcal{R}(S)\subseteq \mathcal{R}(T)$, by Neumann series, we have
	\begin{equation}\label{NeumannS}
	(I+T^\dagger S)^{-1}T^\dagger =\sum_{n=0}^\infty (-T^\dagger S)^{n}T^\dagger =\sum_{n=0}^\infty T^\dagger (-ST^\dagger)^n =T^\dagger (I+ST^\dagger)^{-1}.
	\end{equation}

For the last inequality, we can see that
	\begin{align*}
	(T+S)^\dagger -T^\dagger&=(I+ST^\dagger)^{-1}T^\dagger-(I+ST^\dagger)(I+ST^\dagger)^{-1}T^\dagger \\
	&=[I-(I+ST^\dagger)](I+ST^\dagger)^{-1}T^\dagger\\
	&= (-ST^\dagger)(I+ST^\dagger)^{-1}T^\dagger.
	\end{align*}
	Hence we get the desired inequality.
	
	Now, the fact that $\mathcal{H}=\mathcal{R}(T+S)\oplus\mathcal{N}(T+S)=\mathcal{R}(T)\oplus\mathcal{N}(T)$, and $T$ is invertible from $\mathcal{R}(T)$ to $\mathcal{R}(T)$ , we applied the \cite[Proposition 3.9.]{terElsta2015}, to obtain the last estimate.
\end{proof}
Similarly, we have
\begin{theorem}\label{MPInvsumK}Let $T$ is  m-accretive operator and $S$ is bounded and accretive, then the Theorem \ref{MPInvsum} hold true, if $\mathcal{R}(T)$ is closed, $\mathcal{N}(T)\subseteq \mathcal{N}(S)$ and $\left\| ST^\dagger \right\|  <1$.	
\end{theorem}
\begin{proof} (1)  Since the operator $T+S$ is m-accretive, its adjoint operator $(T+S)^*=T^* +S^*$ is also  m-accretive.
	
(2)	 If $\mathcal{R}(T^*)$ is closed and $\mathcal{N}(T)\subseteq \mathcal{N}(S)$,  then it is obvious that the last inclusion gives  $\mathcal{R}(S^*)\subseteq \mathcal{R}(T^*)$. Also, the condition $\left\| ST^\dagger \right\|  <1$ implies $\left\| (ST^\dagger)^* \right\| =\left\|(T^*)^\dagger S^*  \right\|  <1$, and conversely. Hence by Theorem \ref{MPInvsum}, $\mathcal{R}(T+S)=\mathcal{R}(T^*+S^*)=\mathcal{R}(T^*)=\mathcal{R}(T)$ is closed and $\mathcal{N}(T+S)=\mathcal{N}(T^*+S^*)=\mathcal{N}(T^*)=\mathcal{N}(T)$. Now, we proceed as in the proof of Theorem \ref{MPInvsum}.
\end{proof}
	
\begin{remark}\label{Remark; MPInvsum}Recall that the reduced minimum modulus of a non-zero operator $T$ is defined by
	$$ \gamma (T) = \inf\{ \left\|Tx\right\| : x\in N(T)^{\perp}\cap \mathcal{D}  (T), \left\|x\right\|=1 \}.$$
	If $T = 0$ then we take $\gamma (T)= \infty$. Note that (see \cite{Kato}), $\mathcal{R}  (T)$ is closed if an only if $\gamma (T)>0$. In that case, $\gamma (T)=\dfrac{1}{\left\|T^\dag \right\| }$, where $T^\dag$ is the Moore-Penrose inverse of $T$. Let us remark that if we assume that $\left\| S\right\|  <\dfrac{1}{\gamma (T)}$ instead the condition $\left\| T^\dagger S\right\|  <1$, then the Theorem \ref{MPInvsum}  hold true.
\end{remark}
\begin{proposition}\label{thm:T2macc} Let  $T$ be an accretive such that $T^2$  is  m-accretive. Then
\begin{itemize}
	\item [(i)] $T$ is   m-accretive. Further, if $T$ is $\theta$-accretive with $\theta<\pi/4$, then $T^2$ is m-$2\theta$-accretive.
	\item [(ii)] If $\mathcal{R}(T)$ is closed, then  $\mathcal{R}(T^2)$ is closed and $\gamma(T^2)\geq \dfrac{\gamma(T)^2}{2}$.
\end{itemize}	
	 	
\end{proposition}
\begin{proof} (i) Since $T$ is an accretive operator, by \cite[Theorem 1.2]{Hayashi2017}, we have  , 
	\begin{equation}\label{BBB2}
	\left\|T x\right\|^2 \leq  \nu \left\|x\right\|^2+\dfrac{1}{\nu}\left\|T^2 x\right\|^2,
	\end{equation}
	for all  $x\in  \mathcal{D}  (T^2)$ and an arbitrary $\nu>0$. Choosing $\nu > 0 $ so large that $\dfrac{1}{\nu} < 1$, we obtain $T$ is $T^2$-bounded with lower bound $< 1$. Then $T^2+T$ with domain $\mathcal{D}  (T^2 )$ is 	m-accretive. Now, let us remark that
	$$(\dfrac{1}{4}I+T^2+T)x=(\dfrac{1}{2}I+T)^2x$$
	for all  $x\in  \mathcal{D}  (T^2)$. Since the operator on the left-hand side is invertible, then $(\frac{1}{2}I+T)^2$ is invertible, so  $\frac{1}{2}I+T$ is also invertible. It follows that $T$ is m-accretive. Now, we applied \cite[Theorem 4.]{deLaubenfels1987}.

	(ii) By the Landau-Kolmogorov inequality, \cite[Theorem.]{Kato71}, applied to $T$, we have
	$$\left\|Tx\right\|^2\leq 2 \left\|T^2x\right\|\left\|x\right\|, $$
	for all $x\in\mathcal{N} (T^2)^{\perp}\cap \mathcal{D} (T^2)$. It follows that
	$$ \left\|x\right\|^2 \gamma(T)^2\leq \left\|Tx\right\|^2\leq 2 \left\|T^2x\right\|\left\|x\right\|, $$
	and hence
	$$ \left\|T^2x\right\| \geq \dfrac{\gamma(T)^2}{2}  \left\|x\right\|, $$
	for all $x\in\mathcal{N} (T^2)^{\perp}\cap \mathcal{D} (T^2)$. Now by the definition of $\gamma(T^2)$ we obtain $\gamma(T^2)\geq \dfrac{\gamma(T)^2}{2}$.	
\end{proof}
By Proposition \ref{thm:T2macc}, Theorem \ref{MPInvsum} and Theorem \ref{MPInvsumK},
\begin{corollary}\label{cor:prturb002} Let $T^2$  be m-accretive, $T$ accretive with closed range and $S$  is  bounded and accretive. If    $\mathcal{R}(S)\subseteq \mathcal{R}(T)$ and $\left\| (T^\dagger)^2 S\right\|  <1$. (or $\mathcal{N}(T)\subseteq \mathcal{N}(S)$ and $\left\|  S(T^2)^\dagger\right\|  <1$ ).
	Then \begin{itemize}
		\item $\mathcal{R}(T^2+S)=\mathcal{R}(T^2)=\mathcal{R}(T)$ is closed and $\mathcal{N}(T^2+S)=\mathcal{N}(T^2)=\mathcal{N}(T)$.
		\item  $\mathcal{H}=\mathcal{R}(T^2+S)\oplus\mathcal{N}(T^2+S)$.
		\item $T^2+S$  is an  EP  operator, and
		$$(T^2+S)^\dagger=(I+(T^\dagger)^2 S)^{-1}(T^\dagger)^2 =(T^\dagger)^2
		(I+S(T^\dagger)^2)^{-1}.$$ 
	\end{itemize}
	If further,	$T$ is injective, then $T$, $T^2$ and $T^2+S$  are   invertible, and 
	
	$$ (T^2+S)^{-1} =(I+T^{-2}S)^{-1}T^{-2}=T^{-2}
	(I+ST^{-2})^{-1}.$$ 
	
\end{corollary}
\begin{proof} Since $T^2$  is  m-accretive, then $T$ is also m-accretive by Proposition \ref{thm:T2macc}. Hence $\mathcal{N}(T^2)=\mathcal{N}(T)$ and $\mathcal{R}(T^2)=\mathcal{R}(T)$. By \cite[Lemma 5.5]{Benharrat2020}, we have $(T^\dagger)^2=(T^2)^\dagger$. Now the result is obtained by Theorem \ref{MPInvsum} and Theorem \ref{MPInvsumK}.
\end{proof}	

\begin{remark}
	\begin{itemize}
		\item By Remark \ref{Remark; MPInvsum}, the Corollary \ref{cor:prturb02} hold true if we assume $\left\| T^\dagger S\right\|  <\gamma(T)$ instead $\left\| (T^\dagger)^2 S\right\|  <1$.
		\item \^Ota showed in  \cite[Theorem 2.1]{Ota1984} that, if $T$ is closed and  an accretive such that there is a positive integer $n$ with $\mathcal{D}(T^n)$ is dense in $\mathcal{H}$ and $\mathcal{R}(T^n)\subset \mathcal{D}(T)$, then $T$ is bounded . In particular, for a  closed   and  accretive operator $T$,   if $\mathcal{R}(T)$ is contained in $\mathcal{D}(T)$, or in  $\mathcal{D}(T^*)$, then $T$ is automatically bounded, see also \cite[Theorem 3.3]{Ota1984}.
		\item In general, if $T^2$ is m-accretive; then $T$   fails to be accretive.
		Take $T=i\dfrac{d}{dx}$ on $L^2(\mathbb{R})$. The operator $T$ has its spectrum
		on both sides of the origin. But $T^2=-\dfrac{d^2}{dx^2}$ is a nonegative selfadjoint operator.
	\end{itemize}
	
\end{remark}	


\section{An application to quadratic  operator pencil}

Consider in the Hilbert space $\mathcal{H}$ the  following quadratic  operator pencil
\begin{equation}\label{qop}
\mathcal{Q}(\lambda)=\lambda^2 I-2 \lambda T-S,
\end{equation}
 where $\lambda\in \mathbb{C}$ is the spectral parameter and the two operators  $T$ and $ S $   with domain $\mathcal{D} (T)$ and $\mathcal{D} (S)$, respectively.

One of the approaches to study the spectral properties of  quadratic  operator pencil \eqref{qop} consists of the reducing them to a first order system in a suitable space, \cite{Engstrom2017,Rodman1989}. However, as it was pointed out in \cite{Xiao1998}, this way may be unpractical in the situation when the space is difficult to construct or it is complicated for applications. In addition, as it was mentioned in \cite{XiaoLiang1998,XiaoLiang1999,XiaoLiang2003}, the direct treatment of higher order problems allows to get more general results. An useful approach to address the study and computation of the spectral structure of quadratic  operators pencil is through the use of factorization. This  method  was developed by  Krein and  Langer \cite{KreinL64}, for a   quadratic pencils of self-adjoint operators and by Langer \cite{Langer1976}, for self-adjoint polynomial operator pencils, see also \cite{Langer1976,Langer2004,Harbarth79}. The main idea of this approach consist to   factoring them and studying the spectral properties of the factors. Of particular interest is the separation  of spectral values of $Q_n$ between the  spectra of the  roots. Such separation may be complicated, 
even in the case of eigenvalues, see   \cite{Shkalikov89} and references therein.  The problem is of great importance in spectral theory of such general operators. Moreover, its understanding is crucial in the study of performance properties  of many systems.

The purpose of this section is to extend some earlier factorization  results  essentially
those given in \cite{KreinL64,Harbarth79,Markus88,Moller2015} for the  self-adjoint quadratic operators case, to   \eqref{qop} based on the perturbation theory of  accretive operators together with the uniquely determined  fractional powers of the maximal accretive operators. We also obtain a criterion in order that the linear factors, into which the pencil splits,  generates a holomorphic  semi-group of contraction operators.


We mention that if $T$ is $m$-accretive, then for each $\alpha\!\in\!(0,1)$  the fractional powers $T^\alpha$, $0<\alpha<1$,  are defined by the following Balakrishnan formula, see \cite{Balakrishnan60},
$$
T^\alpha x =\dfrac{\sin (\pi \alpha )}{\pi} \int_{0}^{\infty}\lambda^{\alpha -1}T(\lambda +T)^{ -1}xdt,
$$
for all $x\in \mathcal{D}  (T)$. The operators $T^{\alpha}$ are m-$(\alpha\pi)/2$-accretive  and, if $\alpha\!\in\! (0,1/2)$, then
$\mathcal{D} (T^{\alpha})= \mathcal{D} (T^{*\alpha})$.  It was proved in \cite[Theorem 5.1]{Kato61} that, if $T $ is $m$-accretive, then
$\mathcal{D} (T^{1/2}) \cap \mathcal{D}(T^{*1/2})$ is a core of both $T^{1/2}$  and $ T^{*1/2}$ and the real part ${\rm Re} T^{1/2} :=\!(T^{1/2}\!+\!T^{*1/2})/2 $ defined on $\mathcal{D} (T^{1/2})\cap  \mathcal{D} (T^{*1/2})$ is a selfadjoint operator.
Further,~by~\cite[Corollary ~2]{Kato61},
\begin{equation}
\label{ravhalf}
\mathcal{D}(T)= \mathcal{D}(T^*) \implies \mathcal{D} (T^{1/2})=\mathcal{D} (T^{*1/2})=\mathcal{D} (T^{1/2}_R) = \mathcal{D}[\phi],
\end{equation}
where $\phi$ is the closed form associated with the sectorial operator $T$ via the first representation theorem \cite[Sect.~VI.2.1]{Kato}
and $T_R$ is the non-negative selfadjoint operator associated with the real part of $\phi$ given by
${\rm Re}\, {\phi} := (\phi + \phi^*)/2$, see \cite{Balakrishnan60,Kato,Kato61,Martinez88,Pazy83,Phillips59}. 
 
 Our result  of this section read as follows
\begin{theorem}\label{thm:FAC} Let $T^2$ be m-accretive,  $T$ is accretive and $S$  is accretive bounded operator.Then, we have
	\begin{enumerate}
		\item The operator $\Upsilon  =T^2 +S$ with domain $\mathcal{D} (T^2)$  is  m-accretive.
		\item  $\Upsilon$ admits a fractional powers  $\Upsilon^{\alpha}$  m-$(\alpha \pi/2)$-accretive for each $0<\alpha<1$.
		\item $Z_1=T+\Upsilon^{\frac{1}{2}}$ and $Z_2=T-\Upsilon^{\frac{1}{2}}$ with domain $\mathcal{D} (T)\cap \mathcal{D}(\Upsilon^{\frac{1}{2}})$ are $T^2$-bounded  with lower bound $<1$ and closable operators. Further, $Z_1$ is accretive densely defined.
		
		If further,  $\mathcal{D}(\Upsilon^{\frac{1}{2}})   \subset \mathcal{D} (T) $, then
		
		\item $Z_1$ is m-$\pi/4$-accretive  and for any $\varepsilon>0$, there exists $r>0$, such that $-Z_2-r$ is m-$\pi/4+\varepsilon$-accretive.
		
		In particular, $-Z_1$ and $Z_2+r$   generates  holomorphic  $C_0$-semigroup of contraction operators $\mathcal{T}_1(z)$ and $\mathcal{T}_2(z)$ of angle $\dfrac{\pi}{4}$ and $\pi/4-\varepsilon$, respectively. 
		
		\item The spectra of $Z_1$  and  $Z_2$ are
		separated, $$\sigma(Z_1)\cap \sigma (Z_2)=\emptyset.$$
		
		\item If $T(\mathcal{D} (T^2))\subset \mathcal{D} (T^2)$ and $\Upsilon^{\frac{1}{2}}(\mathcal{D} (T^2))\subset \mathcal{D} (T^2)$, then $Q$ take the following form,
		\begin{equation}\label{z1z2fac}
		\mathcal{Q}(\lambda)x=\dfrac{1}{2}(\lambda I-Z_1)(\lambda I-Z_2)x+\dfrac{1}{2} (\lambda I-Z_2)(\lambda I-Z_1)x,
		\end{equation}
		for all  $x\in \mathcal{D}  (T^2)$.
		
		In particular,  if $TS=ST$ on $\mathcal{D} (T^2)$, then	$Q$ admits the following canonical factorization	
		\begin{equation}\label{z1z2fact}
		\mathcal{Q}(\lambda)x=(\lambda I-Z_1)(\lambda I-Z_2)x=(\lambda I-Z_2)(\lambda I-Z_1)x,
		\end{equation}
		for all  $x\in \mathcal{D}  (T^2)$.
		
	\end{enumerate}
\end{theorem}

\begin{proof} (1) An immediate consequence of Theorem \ref{MPInvsum}-(1).
	
	(2)	$\Upsilon$ admits fractional powers $\Upsilon^{\alpha}$ m-$(\alpha \pi/2)$-accretive for each $0<\alpha<1$,
	see \cite{Balakrishnan60, Kato61}. 
	

	(3) By   (2);  $\Upsilon$ admits unique root $\Upsilon^{\frac{1}{2}}$ m-$(\pi/4)$-accretive operator with $ \mathcal{D}  (T^2)$  is a core of $\Upsilon^{\frac{1}{2}}$. So we  define the following operators
	$$Z_1=T+\Upsilon^{\frac{1}{2}}$$ and $$Z_2=T-\Upsilon^{\frac{1}{2}}$$ with domain $\mathcal{D} (T)\cap \mathcal{D}(\Upsilon^{\frac{1}{2}})$. Both of $Z_1$ and $Z_2$ are densely defined on $\mathcal{H}$ with numerical range is not the whole complex plane, it follows that $Z_1$ and $Z_2$ are closable operators.  Now, we prove that $Z_1$ and $Z_2$  are $T^2$-bounded with lower bound $<1$.  
	
	 By Proposition \ref{thm:T2macc}, $T$ is also m-accretive. Hence, by  \cite[Theorem 1.2]{Hayashi2017}, we have for an arbitrary $\rho_1>0$ and  $\rho_2>0$, 
	\begin{equation}\label{A2A}
	\left\|T x\right\|^2 \leq \rho_1 \left\|x\right\|^2+\dfrac{1}{\rho_1} \left\|T^2x\right\|^2,
	\end{equation}	
	and 
	\begin{equation}\label{lamdas}
	\left\|\Upsilon^{\frac{1}{2}}x\right\|^2\leq  \rho_2\left\|x\right\|^2+\dfrac{1}{\rho_2}\left\|\Upsilon x\right\|^2,
	\end{equation}
	for all $ x\in \mathcal{D}  (T^2)$ (cf. by \cite[Chap. 2, Theorem 6.10]{Pazy83}). 
	
		By  \eqref{A2A} and \eqref{lamdas},  it follows that 
	\begin{align*}
	\left\|Z_ix \right\|^2 & \leq 2 \left\|Tx\right\|^2 + 2\left\|\Lambda^{\frac{1}{2}}x\right\|^2 \\
	& \leq 2 (\rho_1 \left\|x\right\|^2+\dfrac{1}{\rho_1}\left\|T^2 x\right\|^2)+2 (\rho_2\left\|x \right\|^2+\dfrac{1}{\rho_2}\left\|\Lambda x\right\|^2) \\
	& \leq 2 (\rho_1 \left\|x\right\|^2+\dfrac{1}{\rho_1}\left\|T^2 x\right\|^2)+2 \rho_2\left\|x \right\|^2+\dfrac{4}{\rho_2}(\left\|T^2 x\right\|^2 +\left\|Sx\right\|^2) \\
	& \leq  2(\rho_1 +\rho_2 + \dfrac{2\left\|S\right\|^2}{\rho_2}) \left\|x\right\|^2+2 (\dfrac{1}{\rho_1} + \dfrac{2}{\rho_2})\left\|T^2x\right\|^2\\
	&\leq \nu_1 \left\|x\right\|^{2}+ \nu_2\left\|T^2x\right\|^2,
	\end{align*}
	for some   $\nu_1, \nu_2>0$, $i=1,2$ and  all $ x\in \mathcal{D}  (A^2)$. Since $\rho_1$ and $\rho_2$ are arbitrary, we can choose  $\nu_2<1$.

	(4) Now assume  that $\mathcal{D}(\Upsilon^{\frac{1}{2}}) \subset \mathcal{D} (T) $. It follows that
	\begin{equation}\label{AL1sur2}
	\left\|Tx\right\| \leq  a  \left\|x\right\|+ b\left\|\Upsilon^{\frac{1}{2}}x\right\|
	\end{equation}
	for all $x\in \mathcal{D} (\Upsilon^{\frac{1}{2}})$ and for some nonnegative constants $a$ and $b$.	By \eqref{lamdas}, we obtain
	\begin{equation*}
	\left\|Tx\right\|^2 \leq  2a (1+ \rho_2) \left\|x\right\|^2+\dfrac{2b }{\rho_2}\left\|\Upsilon x\right\|^2,
	\end{equation*}
	for all $ x\in \mathcal{D}  (\Upsilon)$ and an arbitrary  $\rho_2>0$. Thus
	\begin{equation*}
	\left\|T (t+ \Upsilon^{\frac{1}{2}})^{-1}x\right\|^2 \leq  2a(1+ \rho_2) \left\|(t+ \Upsilon^{\frac{1}{2}})^{-1}x\right\|^2+\dfrac{2b}{\rho_2}\left\|\Upsilon (t+ \Upsilon^{\frac{1}{2}})^{-1}x\right\|^2,
	\end{equation*}
	for all $x\in \mathcal{H}$. 
	
	Hence
	\begin{equation*}
	\left\|T (t+ \Upsilon^{\frac{1}{2}})^{-1}\right\|^2 \leq  \dfrac{2a}{t^2}(1+  \rho_2) +\dfrac{2b}{\rho_2}\left\|\Upsilon (t+\Upsilon^{\frac{1}{2}})^{-1}\right\|^2.
	\end{equation*}
	Letting $t$ to $+\infty$, we assert that
	$$	M=\sup_{t>0} \left\| T(t+ \Upsilon^{\frac{1}{2}})^{-1} \right\|<\dfrac{2b}{\rho_2^2}. $$
	(cf. \cite[Proposition 2.12]{Yoshikawa72}). Since $\rho_2$ is arbitrary, we  can  choose it such that $\dfrac{2b}{\rho_2^2}<1$. Since $T$ is m-accretive and $\Upsilon^{\frac{1}{2}}$ is m-$(\pi/4)$-accretive, then $Z_1$  is  m-accretive. By \cite[Theorem IX-1.24] {Kato},	the factor $-Z_1$   generates  holomorphic  $C_0$-semigroup $\mathcal{T}_1(z)$  of angle $\dfrac{\pi}{4}$.
	
	On the other hand, since $\Upsilon^{\frac{1}{2}}$ is m-$(\pi/4)$-accretive  and $-T$ satisfy \eqref{AL1sur2}, by \cite[Theorem IX-2.4]{Kato},   for any  $\varepsilon>0$, there exist nonnegative constants $r$ and $s$ such that $a, b<s$ and $(Z_2 +r)$ is the generator of
	holomorphic  $C_0$-semigroup $\mathcal{T}_2(z)$  of angle $\dfrac{\pi}{2}-\dfrac{\pi}{4}-\varepsilon$. This implies that $-(Z_2+r)$ is m-$\psi$-accretive with $\psi=\dfrac{\pi}{4}+\varepsilon$.

	(5) It follows from the item (4),  $W(Z_1)$ is contained in the right half complex plan and $W(Z_2)$ in the let side with a  non zero distance between their closure.

	(6) We have $\mathcal{D} (T^2)\subset \mathcal{D} (Z_1)=\mathcal{D} (Z_2)=\mathcal{D} (\Upsilon^{\frac{1}{2}})\subset \mathcal{D} (T)$. 
	
	The fact that $T(\mathcal{D} (T^2))\subset \mathcal{D} (T^2)$ and $\Upsilon^{\frac{1}{2}}(\mathcal{D} (T^2))\subset \mathcal{D} (T^2)$, we have $\mathcal{D} (T^2)\subset \mathcal{D}(T\Upsilon^{\frac{1}{2}})$, $\mathcal{D} (T^2)\subset \mathcal{D}(\Upsilon^{\frac{1}{2}}T)$ and
	$\mathcal{D} (T^2)\subset \mathcal{D}(Z_1^2)$. Now, by items (1), (2) and (3), we can easily verify  that
	\begin{equation*}
	Z_1^2x -T Z_1x-Z_1Tx-Sx=0,
	\end{equation*}
	for all  $x\in  \mathcal{D} (T^2)$, hence on $\mathcal{D} (T^2)$, we have
	\begin{align*}
	Q(\lambda ) &=Q(\lambda)-(Z_1^2 -T Z_1-Z_1T-S)\\
	&= \lambda^2 I-2 \lambda T-S- Z_1^2 +T Z_1+Z_1T+S\\
	& = \lambda^2 I - Z_1^2 -T(\lambda -Z_1 )- (\lambda -Z_1 )T\\
	&=\dfrac{1}{2}(\lambda -Z_1 )(\lambda I + Z_1-2T)+\dfrac{1}{2}(\lambda I + Z_1-2T)(\lambda -Z_1 )\\
	&=\dfrac{1}{2}(\lambda I-Z_1)(\lambda I-Z_2)+\dfrac{1}{2} (\lambda I-Z_2)(\lambda I-Z_1).
	\end{align*}
	This gives the form \eqref{z1z2fac}. Now, if $TS=ST$ on $\mathcal{D} (T^2)$, then $\Upsilon T=T\Upsilon$. Thus $\Upsilon^{\frac{1}{2}}$ commutes with $T$ on $\mathcal{D} (T^2)$, wich implies that  \eqref{z1z2fact}.

\end{proof}



Now the fact that $\mathcal{N}(\Upsilon^{\frac{1}{2}})=\mathcal{N}(\Upsilon)$ and $\overline{\mathcal{R}(\Upsilon^{\frac{1}{2}})}=\overline{\mathcal{R}(\Upsilon)}$, by Corollary \ref{cor:prturb002} and Theorem \ref{MPInvsumK}, we have
\begin{corollary}\label{cor:prturb02} Let $T^2$  be m-accretive, $T$ accretive with closed range and $S$  is  bounded and accretive. If    $\mathcal{R}(S)\subseteq \mathcal{R}(T)$ and $\left\| (T^\dagger)^2 S\right\|  <1$ (or $\mathcal{N}(T)\subseteq \mathcal{N}(S)$ and $\left\|  S(T^2)^\dagger\right\|  <1$ ).
	Then \begin{itemize}
		\item $\mathcal{R}(Z_1)=\mathcal{R}(Z_2)=\mathcal{R}(\Upsilon^{\frac{1}{2}})=\mathcal{R}(T)$ is closed and $\mathcal{N}(Z_1)=\mathcal{N}(Z_2)=\mathcal{N}(\Upsilon^{\frac{1}{2}})=\mathcal{N}(T)$.
		\item  $\mathcal{H}=\mathcal{R}(\Upsilon^{\frac{1}{2}})\oplus\mathcal{N}(\Upsilon^{\frac{1}{2}})$.
		\item $Z_1$, $Z_2$ and $\Upsilon^{\frac{1}{2}}$  are  EP  operators.
		\item If further,	$T$ is injective, then $Z_1$, $Z_2$ and $\Upsilon^{\frac{1}{2}}$  are   invertible.
	\end{itemize}
\end{corollary}
An immediate consequence of this corollary,   the operator  $Z_1$  and $Z_1$ are written in a matrices form with respect to mutually orthogonal subspaces decomposition as follows
$$Z_1 =\begin{bmatrix}
A& 0 \\
0& 0
\end{bmatrix} : \begin{bmatrix}
\mathcal{R}(T)\\
\mathcal{N}(T)
\end{bmatrix}\longrightarrow  \begin{bmatrix}
\mathcal{R}(T)\\
\mathcal{N}(T)
\end{bmatrix};
$$
and
$$Z_2 =\begin{bmatrix}
B& 0 \\
0& 0
\end{bmatrix} : \begin{bmatrix}
\mathcal{R}(T)\\
\mathcal{N}(T)
\end{bmatrix}\longrightarrow  \begin{bmatrix}
\mathcal{R}(T)\\
\mathcal{N}(T)
\end{bmatrix};
$$
with $A$ and $B$ on  $\mathcal{R}(T) \cap \mathcal{D}(T)$ are injective  operators    with  closed range. In this case, if $AB=BA$ , we have 
$$\mathcal{Q}(\lambda) =\begin{bmatrix}
(\lambda-A)(\lambda-B)& 0 \\
0& \lambda^2 
\end{bmatrix} : \begin{bmatrix}
\mathcal{R}(T)\\
\mathcal{N}(T)
\end{bmatrix}\longrightarrow  \begin{bmatrix}
\mathcal{R}(T)\\
\mathcal{N}(T)
\end{bmatrix}.
$$
Also,  $\mathcal{Q}(0)$ is Moore-penrose invertible and $\mathcal{Q}^{\dagger}(0)=Z_1^{\dagger} Z_2^{\dagger}=A^{-1}B^{-1}$. If $T$ is injective, then $\mathcal{Q}(0)$ is invertible and $\mathcal{Q}^{-1}(0)=Z_1^{-1} Z_2^{-1}$. Consequently, $S$ is also invertible.
The block Vandermonde operator corresponding to $Z_1 , Z_2$ is given by 
$$\mathcal{V}(Z_1 , Z_2)=\begin{bmatrix}
I& I\\
Z_1 & Z_{2}
\end{bmatrix}.$$
We have,

$$\mathcal{V}(Z_1 , Z_2)=\begin{bmatrix}
I& 0\\
Z_2 & I
\end{bmatrix}\begin{bmatrix}
I& 0\\
0 & Z_{1}-Z_2
\end{bmatrix}\begin{bmatrix}
I& I\\
0 & I
\end{bmatrix},$$
where the left and right factors on the right-hand side are invertible. So $\mathcal{V}(Z_1 , Z_2)$ is invertible if and only if $2\Upsilon^{\frac{1}{2}}=(Z_1-Z_2)$ is invertible. Now, we apply \cite[Corollary 29.12, Corollary 29.13 and Remark 29.14] {Markus88} taking in account that $\mathcal{N} (\Lambda^{\frac{1}{2}})= \mathcal{N} ((\Lambda^{\frac{1}{2}})^*)$, we obtain,
$$\sigma(Z_1)\cap \sigma(Z_2)=\emptyset\qquad \text{ and } \qquad \sigma(Z_1)\cup \sigma(Z_2)=\sigma (Q(.)).$$

\section{An application to a second order  linear boundary value problem}%
Denote $[0, +\infty)$ by $\mathbb{R}_+$, and let $C^k(\mathbb{R}_+, \mathcal{D})$ be the
set of all $k$-times (strongly) continuously differentiable functions mapping $\mathbb{R}_+$,
into $\mathcal{D}\subseteq \mathcal{H}$. 
In this section we consider the following  abstract second order  linear boundary value problem,
\begin{align}
u''(t)&-2Tu'(t)-Su(t)=0, \qquad t\in (0, 1),\label{2BVP}\\
u(0)=&u_0, \qquad u(1)=u_1 \label{2BVPC1}
\end{align}
where $u'=\dfrac{d}{dt}$.
\begin{theorem}\label{thm:sol2BVP}  Let $T^2$  be m-accretive, $T$ is accretive  and $S$  is  bounded and accretive.   	Assume that 
	\begin{enumerate}
		\item $\mathcal{D}(\Upsilon^{\frac{1}{2}})   \subset \mathcal{D} (T) $.
		\item $T(\mathcal{D} (T^2))\subset \mathcal{D} (T^2)$ and $\Upsilon^{\frac{1}{2}}(\mathcal{D} (T^2))\subset \mathcal{D} (T^2)$.
		\item $T$ commutes with $\Upsilon^{1/2}$ on $\mathcal{D}(T^2)$.
	\end{enumerate}
	Then of any   constant vectors
	$u_0$ ,$u_1 \in\mathcal{D} (T)$  the vector valued function,
	$$u(t)=e^{-(1-t) Z_1}x_0+e^{tZ_2}x_1, \qquad t\in (0, 1),$$
	with
	$$x_0=(I-e^{-2 \Upsilon^{\frac{1}{2}}})^{-1}\left[-e^{Z_2} u_0+ u_1 \right] $$
	and
	$$x_1=(I-e^{-2 \Upsilon^{\frac{1}{2}}})^{-1}\left[ u_0-e^{-Z_1}u_1 \right] $$
	is uniquely determined solution of \eqref{2BVP}-\eqref{2BVPC1}, with $u(.)\in C^{\infty}((0, 1), \mathcal{H}) \cap C^1((0, 1), \mathcal{D} (T))$.
\end{theorem}
\begin{proof}  Under the assumptions, by Theorem \ref{thm:FAC}, the factors  $-Z_1$  and $Z_2$ 
	generates bounded  holomorphic  $C_0$-semigroups. By \cite[Lemma 2.38]{Davies1980},
	$$x(t)=e^{-(1-t)Z_1}x_0\in \mathcal{D}  (Z_1)^k)$$
	and
	$$y(t)=e^{tZ_2}x_1\in \mathcal{D}  (Z_2)^k) ,$$
	for all  $k\in\mathbb{N}$, $x_0, x_1\in \mathcal{H}$ and 	$t\in (0, 1)$. This implies that
	$$u(t)=x(t)+y(t)\in \mathcal{H},$$
	$$u'(t)=Z_1x(t)+Z_2 y(t)\in \mathcal{D} (T) $$
	and
	$$u^{(2)}(t)=Z_1^2x(t)+Z_2^2 y(t)$$
	for all $t\in (0, 1)$. We can easily see  that $u$ verifies \eqref{2BVP}. Since $(I-e^{-2 \Upsilon^{\frac{1}{2}}})^{-1}$ and $(I-e^{-2 \Upsilon^{\frac{1}{2}}})^{-1}$ exist and bounded  on $\mathcal{H}$. Thus,  we have    
		$$u(0)=e^{- Z_1}x_0+x_1=u_0, $$ and
	$$u(1)=x_0+e^{-Z_2}x_1=u_1, $$
	This completes the proof.
\end{proof}
\begin{example}
	Let $\Omega $ is a smooth bounded domain of $\mathbb{R}^n$, $\Gamma$ be the
	boundary of $\Omega $, and $\xi\in \mathbb{C}$ with ${\rm Re}(\xi)\geq 0$, $\eta\geq 0$, $\eta_1\in \mathbb{R}$.  We consider the  following initial-boundary value problem in $L^2(\Omega)$,
	\begin{align}
	u'' (t, x)& -2 (\eta \Delta - i \eta_1 \Delta^2 )  u'(t, x)-\xi u(t, x)=0, \qquad (t, x) \in (0,1) \times \Omega ,\label{mIBVP}\\
	u(0, x)=&u_0, \qquad u(1, x)=u_1, \qquad x \in \Omega , \label{mIBVPC1}\\
	u_{|\Gamma}=&\Delta u_{|\Gamma}=0, \label{mIBVPC2}
	\end{align}
	where $\Delta$ denotes the Laplacian. It is known that $-\Delta$ with domain $H^{1}_0(\Omega)\cap H^{2}(\Omega)=\{ u\in H^2(\Omega) ; u_{|\Gamma}=0\}$ is a positive, self-adjoint operator on $L^2(\Omega)$. Its $L^2(\Omega)$ - normalized eigenfunctions are denoted $w_j$, and its eigenvalues counted with their multiplicities are denoted $\lambda_j$: 
	\begin{equation}
	-\Delta w_j = \lambda_j w_j.
	\label{ef}
	\end{equation}
	It is well known that $0<\lambda_1\leq...,\leq \lambda_j \longrightarrow \infty$. Functional calculus can be defined using the eigenfunction expansion. In particular; if  we denote by $(-\Delta )^{\alpha}$ the fractional powers of the Dirichlet Laplacian, with $0\le \alpha\le 1$, then
	\begin{equation*}
	\left(-\Delta\right)^{\alpha}u = \sum_{j=1}^{\infty}\lambda_j^{\alpha} u_j w_j
	\label{funct}
	\end{equation*}
	with 
	\[
	u_j =\int_{\Omega}u(x)w_j(y)dx
	\]
	for $u\in{\mathcal{D}}\left( (-\Delta )^{\alpha} \right)= \{u :  \; (\lambda_j^{\alpha}u_j)\in \ell^2(\mathbb {N})\}$. Also, if $m\alpha$ is an  integer, then
	\begin{equation}\label{domain}
	\mathcal{D}\left( (-\Delta )^{m\alpha} \right)=H^{m\alpha}_0(\Omega)\cap H^{2m\alpha}(\Omega),\qquad m\geq 1.
	\end{equation}
	See \cite{Lions72}, for more details.
	
	Set	$$ T= \eta \Delta -i \eta_1 \Delta^2,\qquad  \mathcal{D} (T)= \mathcal{D} (\Delta^2)=\{ u\in H^4(\Omega) ; u_{|\Gamma}=\Delta u_{|\Gamma}=0\}.$$

	and
		$$Su= \xi  u, \qquad \mathcal{D} (S)= L^{2}(\Omega) $$ 
	Then the abstract version of problem \eqref{mIBVP}-\eqref{mIBVPC1} takes the form \eqref{2BVP}-\eqref{2BVPC1}. We have,
	\begin{itemize}
		\item $S$ is bounded and accretive. 
		\item $T$ is m-accretive.
		\item The operator 
		$$T^2=\eta^2 \Delta^2 + \eta_1^2 \Delta^4 -2i\eta \eta_1 \Delta^3 $$   with $$\mathcal{D} (T^2)=\mathcal{D} (\Delta^4)=\{ u\in H^8(\Omega) ; u_{|\Gamma}=\Delta u_{|\Gamma} =\Delta^2 u_{|\Gamma}=\Delta^3 u_{|\Gamma}=0\}$$  is m-accretive. Indeed, $T^2$ is densely defined and closed, and $(T^2)^*=\eta^2 \Delta^2 + \eta_1^2 \Delta^4 +2i\eta \eta_1 \Delta^3$ is accretive.
		\item The operator $$\Upsilon=T^2 +S= \eta^2 \Delta^2 + \eta_1^2 \Delta^4 -2i\eta \eta_1 \Delta^3 +\xi,$$ with domain $\mathcal{D} (\Delta^4)$, is m-accretive. 
		\item $\Upsilon$ admits a square root  $\Upsilon^{1/2}$  m-$(\pi/4)$-accretive with $\mathcal{D} (\Delta^4)$ is a core of $\Upsilon^{1/2}$.
		\item The operators factors $$Z_1=T+\Upsilon^{\frac{1}{2}}=\eta \Delta -i \eta_1 \Delta^2+(\eta^2 \Delta^2 + \eta_1^2 \Delta^4 -2i\eta \eta_1 \Delta^3 +\xi)^{\frac{1}{2}} $$ 
		and 
		$$Z_2=T-\Upsilon^{\frac{1}{2}}=\eta \Delta -i \eta_1 \Delta^2 -(\eta^2 \Delta^2 + \eta_1^2 \Delta^4 -2i\eta \eta_1 \Delta^3 +\xi)^{\frac{1}{2}}$$ with domain $\mathcal{D} (\Delta^2)$ are closed operators.
		\item By an argument of functional calculus, we  obtain
		$$Tu=  -\sum_{j=1}^{\infty}(\eta + i\eta_1 \lambda_j)\lambda_j u_j w_j $$
		for $u\in \mathcal{D} (\Delta^2)$.
		If  $(\eta,  \eta_1)\neq (0,0)$ then $T$ is injective. Thus, is invertible with
		$$T^{-1}u=  -\sum_{j=1}^{\infty} \dfrac{1}{(\eta + i\eta_1 \lambda_j)\lambda_j} u_j w_j $$
		for $u\in \mathcal{D} (\Delta^2)$.
		\item Assume that $(\eta,  \eta_1)\neq (0,0)$ and $\xi\neq 0$. If  $ \sum_{j=1}^{\infty} \dfrac{1}{(\eta + i\eta_1 \lambda_j)^2\lambda_j^2}<\dfrac{1}{\left|\xi \right| }$ then $\left\|  ST^{-2}\right\|  <1$. In fact,
		\begin{align*}
		\left\|  T^{-1}u\right\|^2&=\left\langle T^{-2}u, u \right\rangle \\
		& =\sum_{j=1}^{\infty} \dfrac{1}{(\eta + i\eta_1 \lambda_j)^2\lambda_j^2}\left|\left\langle u, w_j \right\rangle  \right|^2 \\
		&\leq \sum_{j=1}^{\infty} \dfrac{1}{(\eta + i\eta_1 \lambda_j)^2\lambda_j^2} \left\| u\right\|^2\\
		&<\dfrac{1}{\left|\xi \right| } \left\| u\right\|^2.
		\end{align*}
		This implies that $$\left\|  ST^{-2}\right\|\leq\left\|  S\right\|\left\|  T^{-1}\right\|^2= \left|\xi \right|\left\|  T^{-1}\right\|^2<1.$$
		\item If we assume that $(\eta,  \eta_1)\neq (0,0)$,  $\xi\neq 0$ and   $ \sum_{j=1}^{\infty} \dfrac{1}{(\eta + i\eta_1 \lambda_j)^2\lambda_j^2}<\dfrac{1}{\left|\xi \right| }$; then by Theorem
		\ref{MPInvsumK} and Corollary \ref{cor:prturb02}, we conclude that  $\Upsilon^{\frac{1}{2}}$, $Z_1$ and $-Z_2$ are m-$(\pi/4)$-accretive invertible  operators.  In particular, $Z_1$ and $-Z_2$   generates  holomorphic  $C_0$-semigroup of contraction operators $\mathcal{T}_1(z)$ and $\mathcal{T}_2(z)$ of angle $\dfrac{\pi}{4}$. If we assume further $T$ commutes with $\Upsilon^{\frac{1}{2}}$ on $\mathcal{D} (T^2)$. Then all the statements of Theorem \ref{thm:sol2BVP} hold. Consequently, for any pair of  vectors
		$u_0, u_1\in  \mathcal{D} (T)$  the vector valued function,
		$$u(t, x)=e^{-(1-t) Z_1}v_0(x)+e^{tZ_2}v_1(x), \qquad t\in (0, 1),  \quad x\in \Omega$$
		with
		$$v_0(x)=(I-e^{-2 \Upsilon^{\frac{1}{2}}})^{-1}\left[-e^{Z_2} u_0(x)+ u_1(x) \right] \quad x\in \Omega $$
		and
		$$v_1(x)=(I-e^{-2 \Upsilon^{\frac{1}{2}}})^{-1}\left[ u_0(x)-e^{-Z_1}u_1(x) \right] \quad x\in \Omega$$
		is the  unique solution of \eqref{mIBVP}-\eqref{mIBVPC2}.
	\end{itemize}

\end{example}



\end{document}